\magnification=\magstep1
\input amstex
\documentstyle{amsppt}
\hoffset=.25truein
\hsize=6truein
\vsize=8.75truein

\topmatter
\title
On exceptional eigenvalues of the Laplacian for $\Gamma_0(N)$
\endtitle
\abstract
   An explicit Dirichlet series is obtained, which represents
an analytic function of $s$ in the half-plane $\Re s>1/2$ except
for having simple poles  at points $s_j$ that
correspond to exceptional eigenvalues $\lambda_j$
of the non-Euclidean Laplacian for Hecke congruence subgroups
$\Gamma_0(N)$ by the relation $\lambda_j=s_j(1-s_j)$ for
$j=1,2,\cdots, S$.  Coefficients of the Dirichlet series involve
all class numbers $h_d$ of real quadratic number fields.  But,
only the terms with $h_d\gg d^{1/2-\epsilon}$ for sufficiently
large discriminants $d$ contribute to the residues $m_j/2$ of the Dirichlet
series at the poles $s_j$, where $m_j$ is the multiplicity of the
eigenvalue $\lambda_j$ for $j=1,2,\cdots, S$.   This may indicate
(I'm not able to prove yet) that the multiplicity of exceptional
eigenvalues can be arbitrarily large.  On the other hand, by density
theorem [3] the multiplicity of exceptional eigenvalues is bounded
above by a constant depending only on $N$.
\endabstract
\author
 Xian-Jin Li
\endauthor
\address
Department of Mathematics, Brigham Young University, Provo, Utah
84602 USA
\endaddress
\email
xianjin\@math.byu.edu
\endemail
\thanks
Research supported by National Security Agency H98230-06-1-0061
\endthanks
\endtopmatter
\document

\heading 1.   Introduction
\endheading

   Let $N$ be a positive integer.  Denote by
$\Gamma_0(N)$ the Hecke congruence subgroup of level $N$. The
non-Euclidean Laplacian $\Delta$ on the upper half-plane $\Cal H$ is
given by
$$\Delta= -y^2\left(\frac{\partial ^2}{\partial x^2}+\frac{\partial
^2}{\partial y^2} \right).$$ Let $D$ be the fundamental domain of
$\Gamma_0( N)$.   Eigenfunctions of the discrete spectrum of
$\Delta$ are nonzero real-analytic solutions of the equation
$\Delta\psi=\lambda\psi$ such that $\psi(\gamma z)=\psi(z)$ for all
$\gamma$ in $\Gamma_0( N)$ and such that $\psi$ is square integrable
on $D$ with respect to the Poincar\'e measure $dz$ of the upper
half-plane.

   The Hecke operators $T_n$, $n=1, 2, \cdots$, $(n, N)=1$, which act
in the space of automorphic functions with respect to $\Gamma_0(
N)$, are defined by
$$\left(T_nf\right)(z)=\frac 1{\sqrt n}\sum_{ad=n, \,
0\leq b<d}f\left(\frac {az+b}d\right). $$ It is well-known (see
Iwaniec [3]) that there exists a maximal orthonormal system of
eigenfunctions of $\Delta$ such that each of them is an
eigenfunction of all the Hecke operators. Let $\lambda_j$, $j=1, 2,
\cdots$, be an enumeration in increasing order of all positive
discrete eigenvalues of $\Delta$ for $\Gamma_0(N)$ with an
eigenvalue of multiplicity $m$ appearing $m$ times, and let
$\kappa_j=\sqrt{\lambda_j-1/4}$ with $\Im \kappa_j>0$ if
$\lambda_j<1/4$.

If $\lambda$ is a positive discrete eigenvalue
less than $1/4$, we call it an exceptional eigenvalue.
Let $\lambda_1, \cdots, \lambda_S$ be exceptional eigenvalues
  of the Laplacian $\Delta$ for $\Gamma_0(N)$.

  In 1965, A. Selberg [10] made the following fundamental conjecture.

  \proclaim{Selberg's eigenvalue conjecture}  If $\lambda$ is a
nonzero discrete eigenvalue of the non-Euclidean Laplacian for any
congruence subgroup, then $\lambda \geq 1/4$. \endproclaim

  A. Selberg [10] proved that $\lambda\geqslant 3/16$.
 The best available lower bound
$\lambda\geq 975/4096$ is due to Kim and Sarnak [4].  It was
obtained by combining automorphic lifts sym$^3: GL(2)\to GL(4)$
[5] and sym$^4: GL(2)\to GL(5)$ [4] with families of $L$-functions [8].
We note that if the general functorial conjectures concerning
the automorphic lifts sym$^k: GL(2)\to GL(k+1)$ are true for
all $k>1$, then Selberg's eigenvalue conjecture would follow.

  In this paper, we indicate an elementary approach towards
the Selberg eigenvalue conjecture.  Namely, we prove the following
theorem.

\proclaim{Theorem 1}  Let
$$L(s)=\sum_{k|N}\sum_{d\in \Omega,\, k|u_d}
\prod_{p^{2l}|(d, N/k)}p^l \prod_{p| N/k}(1+({d\over p}))
{h_{dk^2}\ln \epsilon_{dk^2}\over (du_d^2)^s} -\sum_{d\in\Omega}
\frac{h_d\ln \epsilon_d}{(d u_d^2)^s} $$
where $(v_d, u_d)$ is
the smallest positive solution of Pell's equation $v^2-du^2=4$ and
the product on $p^{2l}$ is over all distinct primes $p$ with
$p^{2l}$ being the greatest even $p$-power factor of $(d,N/k)$.
Then $L(s)$ represents an analytic function of $s$ in the
half-plane $\Re s > 1/2$ except for having simple poles at
$s_j=\frac 12- i\kappa_j$, $j=1,2,\cdots, S$.  Moreover, we have
$$m_j=2\,\text{Res}_{s=s_j} L(s)$$
for $j=1,2,\cdots, S$.
\endproclaim

 \proclaim{Corollary 2}  If $N$ is square free, then the series
$$L_1(s)=\sum_{\underset{(m, k)\neq (1,1)}\to {m|N, \,k|N}}
k{\mu((m, k))\over (m,k)} \sum_{d\in \Omega,\,
k|u_d} \left({d\over m}\right) {h_d\ln \epsilon_d\over (du_d^2)^s}$$
represents an analytic function of $s$ in the
half-plane $\Re s > 1/2$ except for having simple poles at
 $s_j=\frac 12- i\kappa_j$, $j=1,2,\cdots, S$,
where $(v_d, u_d)$ is the smallest positive solution of Pell's
equation $v^2-du^2=4$.  Moreover, we have
$$m_j=2\,\text{Res}_{s=s_j} L_1(s).$$
\endproclaim

This work was initiated while the author attended the
workshop on Eisenstein Series and Applications at
the American Institute of Mathematics (AIM), August 15-19, 2005.
The author wants to thank AIM for the invitation of attending
the workshop.

\heading
2.  Proofs of Theorem 1 and Corollary 2
\endheading

 We denote by $h_d$ the class number of
indefinite rational quadratic forms with discriminant d.  Let
$$\epsilon_d=\frac{v_0+u_0\sqrt d}2, $$
where the pair $(v_0, u_0)$ is the smallest positive solution of Pell's
equation $v^2-du^2=4$.  Let $\Omega$ be the set of all the
positive integers $d$ such that $d\equiv 0$ or 1 (mod 4) and such
that $d$ is not a square of an integer.

\proclaim{Lemma 2.1}  Let $d$ and $d_1$ be integers in $\Omega$.
If $d_1=dl^2$, then
$$h_{d_1}\ln \epsilon_{d_1}=l\prod_{p|l}
\left(1-\left({d\over p}\right)p^{-1}\right)h_d\ln \epsilon_d.$$
\endproclaim

\demo{Proof}   The stated identity follows from Dirichlet's
class number formula (see, \S100 of Dirichlet [2])
$$h_{d_1}\ln\epsilon_{d_1}=\sqrt{d_1}L(1, \chi_{d_1})$$
and the identity
$$L(1, \chi_{d_1})=L(1, \chi_d)\prod_{p|l}
\left(1-\left({d\over p}\right)p^{-1}\right). \qed$$

  \enddemo

\proclaim{Lemma 2.2}  Let $d$ and $d_1$ be integers in $\Omega$, and
let $d_1=dl^2$.  Then $\epsilon_{d_1}=\epsilon_d^{\nu_l}$
for a positive integer $\nu_l$.
\endproclaim

\demo{Proof}  If $(v_1, u_1)$ is the smallest positive solution of
Pell's equation
$$v^2-dl^2u^2=4, \tag 2.1 $$
 then
$$\epsilon_{d_1}={v_1+\sqrt{d_1}u_1\over 2}.$$
Let $(v_0, u_0)$ be the smallest positive solution of
Pell's equation
$$v^2-du^2=4. \tag 2.2 $$
By \S85 of Dirichlet [2], all positive solutions $(v, u)$ of (2.2) are
given by the formula
$${v+\sqrt du\over 2}=\left({v_0+\sqrt du_0\over 2}\right)^n$$
for positive integers $n$.  Since $(v_1, lu_1)$ is a positive
solution of (2.2),  there exists a positive integer $\nu_l$ such that
$${v_1+\sqrt {d_1}u_1\over 2}=\left({v_0+\sqrt du_0\over 2}\right)^{\nu_l}.
\qed $$
\enddemo

    We denote the multiplicity of the eigenvalue $\lambda_j$
  by $m_j$ for $j=1,2,\cdots$.

\proclaim{Lemma 2.3}   Let $N$ be any positive
integer, and let
$$L_N(s)=\sum_{k|N}k^{1-2s}\sum_{d\in \Omega}\sum_u
\prod_{p^{2l}|(d, N/k)}p^l
 \prod_{p| N/ k}(1+({d\over p}))
\prod_{p|k}(1-({d\over p})p^{-1})
{h_d\ln \epsilon_d\over (du^2)^s} $$
for $\Re s>1$, where the sum on $u$ is  over all positive integers
$u$ such that $\sqrt{4+d k^2 u^2}\in\Bbb Z$ and where
the product on $p^{2l}$ is over all
distinct primes $p$ with $p^{2l}$ being the greatest even
$p$-power factor of $(d,N/k)$.  Then $L_N(s)$ is analytic for
$\Re s>1$ and has analytic continuation to the
half-plane $\Re s > 1/2$ except for having simple poles
at $s=1$ and  $s_j=\frac 12- i\kappa_j$,
$j=1,2,\cdots, S$.  Moreover, we have
$$m_j=2\text{Res}_{s=s_j} L_N(s)$$
for $j=1,2,\cdots, S$.
\endproclaim

\demo{Proof}  Let
$$h(r)=4^s\sqrt\pi {\Gamma(s-1/2)\over\Gamma(s)}
\int_0^\infty \left(u+{1\over u}+2\right)^{1/2-s}u^{ir-1}du$$
for $\Re s>1/2$.
Then the lemma follows from Theorem 4.3, the proof of
Lemma 5.3, the proof of Theorem 1 in Li [7], and the Selberg trace formula
$$\aligned h(-i/2)&+\sum_{j=1}^\infty h(\kappa_j)m_j\\
&=4^{1/2+s}\sqrt \pi {\Gamma(s-1/2)\over\Gamma(s)}L_N(s) +f(s)
\endaligned $$
for $\Re s>1$ where $f(s)$ is a certain analytic function of $s$
in the half-plane $\Re s\geq 1/2$ except for a possible pole at $s=1/2$
 (see (4.4) of Li [7]).
\qed\enddemo

 \remark{ Remark 2.4}  Siegel [11] proved that
$$\lim_{d\to\infty} \frac{\ln(h_d\ln\epsilon_d)}{\ln d}=\frac 12.
\tag 2.3$$
 \endremark

\proclaim{Lemma 2.5}  Let
$$\bar L_N(s)=\sum_{k|N}\sum_{d\in \Omega,\, k|u_d}
\prod_{p^{2l}|(d, N/k)}p^l \prod_{p| N/k}(1+({d\over p}))
{h_{dk^2}\ln \epsilon_{dk^2}\over (du_d^2)^s} $$
for $\Re s>1$, where $(v_d, u_d)$ is the smallest positive solution
of Pell's equation $v^2-du^2=4$ and
the product on $p^{2l}$ is over all
distinct primes $p$ with $p^{2l}$ being the greatest even
$p$-power factor of $(d,N/k)$.  Then $\bar L_N(s)$ is analytic for
$\Re s>1$ and has analytic continuation to the
half-plane $\Re s > 1/2$ except for having simple poles
at $s=1$ and $s_j=\frac 12- i\kappa_j$,
$j=1,2,\cdots, S$.  Moreover, we have
$$m_j=2\text{Res}_{s=s_j} \bar L_N(s).$$
\endproclaim

\demo{Proof}   By Lemma 2.3, the function
$$L_N(s)=\sum_{k|N}k^{1-2s}\sum_{d\in \Omega}\sum_u
\prod_{p^{2l}|(d, N/k)}p^l
 \prod_{p| N/ k}(1+({d\over p}))
\prod_{p|k}(1-({d\over p})p^{-1})
{h_d\ln \epsilon_d\over (du^2)^s} \tag 2.4$$
 has analytic continuation to the
half-plane $\Re s > 1/2$ except for having simple poles
at $s=1$ and $s_j=\frac 12- i\kappa_j$,
$j=1,2,\cdots, S$,  where the sum on $u$ is  over all positive
solutions of  Pell's equation
$$v^2-d k^2 u^2=4. \tag 2.5 $$

   Let $(v_k, u_k)$ be the smallest positive solution of (2.5).
By \S85 of Dirichlet [2], all positive solutions $(v, u)$ of (2.5) are
given by the formula
$${v+\sqrt dku\over 2}=\left({v_k+\sqrt dku_k\over 2}\right)^n$$
for $n=1,2,\cdots$.  Hence, we have
$$\aligned \sqrt d ku&=\left({v_k+\sqrt dku_k\over 2}\right)^n
\left(1-\left({v_k+\sqrt dku_k\over 2}\right)^{-2n}\right)\\
&>\left({v_k+\sqrt d ku_k\over 2}\right)^n(1+2/\sqrt dk)^{-1}.
\endaligned \tag 2.6$$

  Let $\sigma=\Re s>1/2$, and let $\tau(n)$ be the number of
  positive divisors of an integer $n$.  By (2.6) and (2.3), we have
    $$\aligned &|\sum_{k|N}k^{1-2s}\sum_{d\in \Omega}\sum_{u\neq u_k}
\prod_{p^{2l}|(d, N/k)}p^l
 \prod_{p| N/ k}(1+({d\over p}))
\prod_{p|k}(1-({d\over p})p^{-1})
{h_d\ln \epsilon_d\over (du^2)^s}|\\
&\leq \sum_{k|N}\sqrt{kN}2^{\tau(N)}
\sum_{d\in\Omega}(1+2/\sqrt dk)^{2\sigma}h_d\ln\epsilon_d
\sum_{n=2}^\infty \left({v_k+\sqrt dku_k\over 2}\right)^{-2n\sigma}\\
&\leq \sum_{k|N}\sqrt{kN}2^{\tau(N)}3^{2\sigma+1}
\sum_{d\in\Omega}h_d\ln\epsilon_d
\left({v_k+\sqrt dku_k\over 2}\right)^{-4\sigma}\\
&\leq 16^\sigma N2^{\tau(N)}3^{2\sigma+1}
\sum_{k|N}\sum_{d\in\Omega}d^{1/2+\epsilon-2\sigma}(ku_k)^{-4\sigma}.
  \endaligned $$
Note that $\tau(n)=n^\epsilon$ as $n\to\infty$.  Since, for a fixed positive
integer $v$, there are at most $\tau(v^2-4)$ number of $d$'s in $\Omega$
such that $v^2-du^2=4$ for positive integers $u$, we have
  $$\sum_{d\in\Omega}d^{1/2+\epsilon-2\sigma}(ku_k)^{-4\sigma}
  \leq \sum_{d\in\Omega} (v_k^2-4)^{-\sigma}
  \leq \sum_{v=3}^\infty {\tau(v^2-4)\over (v^2-4)^\sigma}
  \ll \sum_{v=3}^\infty {1\over (v^2-4)^{\sigma-\epsilon}}<\infty
   $$
  for $\sigma>1/2$.  Hence, the series
  $$\sum_{k|N}k^{1-2s}\sum_{d\in \Omega}\sum_{u\neq u_k}
\prod_{p^{2l}|(d, N/k)}p^l
 \prod_{p| N/ k}(1+({d\over p}))
\prod_{p|k}(1-({d\over p})p^{-1})
{h_d\ln \epsilon_d\over (du^2)^s}$$
represents an analytic function of $s$
 in the half-plane $\Re s>1/2$.  It follows from
(2.4) that the function
$$\sum_{k|N}k^{1-2s}\sum_{d\in \Omega}
\prod_{p^{2l}|(d, N/k)}p^l
 \prod_{p| N/ k}(1+({d\over p}))
\prod_{p|k}(1-({d\over p})p^{-1})
{h_d\ln \epsilon_d\over (du_k^2)^s} \tag 2.7$$
has analytic continuation to the
half-plane $\Re s > 1/2$ except for having simple poles
at $s=1$ and $s_j=\frac 12- i\kappa_j$,
$j=1,2,\cdots, S$.

  Next, let $(v_0, u_0)$ be the smallest positive solution of
Pell's equation
$$v^2-du^2=4. \tag 2.8$$
Let $k$ be a divisor of $N$.  If $(v_k, ku_k)$
is a solution of (2.8) different from $(v_0, u_0)$,  then by Lemma 2.2
there exists an integer $n\geq 2$ such that
$${v_k+\sqrt dku_k\over 2}=\left({v_0+\sqrt du_0\over 2}\right)^n.$$
  Hence, we have
$$\aligned \sqrt d ku_k=&\left({v_0+\sqrt du_0\over 2}\right)^n
\left(1-\left({v_0+\sqrt du_0\over 2}\right)^{-2n}\right)\\
&>\left({v_0+\sqrt d u_0\over 2}\right)^n(1+2/\sqrt d)^{-1}.
\endaligned\tag 2.9$$

    By (2.9) and (2.3), we have
    $$\aligned &|\sum_{k|N}k^{1-2s}\sum_{d\in \Omega, \, ku_k\neq u_0}
\prod_{p^{2l}|(d, N/k)}p^l
 \prod_{p| N/ k}(1+({d\over p}))
\prod_{p|k}(1-({d\over p})p^{-1})
{h_d\ln \epsilon_d\over (du_k^2)^s}|\\
&\leq \sum_{k|N}\sqrt{kN}2^{\tau(N)}
\sum_{d\in\Omega,\, ku_k\neq u_0}(1+2/\sqrt d)^{2\sigma}h_d\ln\epsilon_d
 \left({v_0+\sqrt du_0\over 2}\right)^{-2n\sigma}\\
&\leq \sum_{k|N}\sqrt{kN}2^{\tau(N)}9^{\sigma}
\sum_{d\in\Omega, \, ku_k\neq u_0}h_d\ln\epsilon_d
\left({v_0+\sqrt du_0\over 2}\right)^{-4\sigma}\\
&\leq 16^\sigma \tau(N)N2^{\tau(N)}9^{\sigma}
\sum_{d\in\Omega}d^{1/2+\epsilon-2\sigma}u_0^{-4\sigma}.
  \endaligned $$

  Since
  $$\sum_{d\in\Omega}d^{1/2+\epsilon-2\sigma}u_0^{-4\sigma}
  \leq \sum_{d\in\Omega} (v_0^2-4)^{-\sigma}
  \leq \sum_{v=3}^\infty {\tau(v^2-4)\over (v^2-4)^\sigma}
  \leq \sum_{v=3}^\infty {1\over (v^2-4)^{\sigma-\epsilon}}<\infty
   $$
  for $\sigma>1/2$, the series
  $$\sum_{k|N}k^{1-2s}\sum_{d\in \Omega,\, ku_k\neq u_0}
\prod_{p^{2l}|(d, N/k)}p^l
 \prod_{p| N/ k}(1+({d\over p}))
\prod_{p|k}(1-({d\over p})p^{-1})
{h_d\ln \epsilon_d\over (du_k^2)^s}$$
represents an analytic function of $s$
 in the half-plane $\Re s>1/2$.  It follows from
 (2.7) that the function
$$\sum_{k|N}k\sum_{d\in \Omega,\, k|u_0}
\prod_{p^{2l}|(d, N/k)}p^l
 \prod_{p| N/k}(1+({d\over p}))\prod_{p|k}(1-({d\over p})p^{-1})
{h_d\ln \epsilon_d\over (du_0^2)^s} \tag 2.10 $$
has analytic continuation to the
half-plane $\Re s > 1/2$ except for having simple poles
at $s=1$ and  $s_j=\frac 12- i\kappa_j$,  $j=1,2,\cdots, S$,
where the product on $p^{2l}$ is over all
distinct primes $p$ with $p^{2l}$ being the greatest even
$p$-power factor of $(d,N/k)$.
By Lemma 2.1 we can write (2.10) as
$$\sum_{k|N}\sum_{d\in \Omega,\, k|u_0}
\prod_{p^{2l}|(d, N/k)}p^l \prod_{p| N/k}(1+({d\over p}))
{h_{dk^2}\ln \epsilon_{dk^2}\over (du_0^2)^s}. $$

This completes the proof of the lemma. \qed\enddemo

 \demo{Proof of Theorem 1}  It is proved in [6] that the series
$$F(s)=\sum_{d\in\Omega} \frac{h_d\ln \epsilon_d}{d^s}
\sum_{\underset {v^2-du^2=4}\to {u>0}} {1\over u^{2s}},$$
represents an analytic function of $s$ in the half-plane
$\Re\, s> 1/2$ except for having a simple poles at $s=1$.
  By (2.3), (3.4), (3.5), Lemma 3.5,
Lemma 4.1, and Lemma 4.2 of [6], we have that
$$F(s)-h(-i/2) \tag 2.11$$
 is analytic in the half-plane $\Re\, s>1/2$.
 Let $(v_d, u_d)$ be the smallest positive solution of
 Pell's equation $v^2-du^2=4$.  If $u\neq u_d$, then
 $${v+\sqrt d  u\over 2}=\left({v_d+\sqrt d  u_d\over 2}\right)^\nu$$
 for some positive integer $\nu\geq 2$.  Similarly as in (2.9), we
 can obtain that
 $$\sqrt du\geq {1\over 3}\epsilon_d^\nu.$$
  It follows that
 $$\aligned \left|\sum_{d\in\Omega} \frac{h_d\ln \epsilon_d}{d^s}
\sum_{\underset {v^2-du^2=4}\to {u\neq u_d}} {1\over u^{2s}}\right|
&\leq 9\sum_{d\in\Omega} \frac{h_d\ln \epsilon_d}{\epsilon_d^{4\sigma}}\\
&\leq  2^{4\sigma}9 \sum_{d\in\Omega}d^{1/2+\epsilon-2\sigma} u_d^{-4\sigma}
<\infty \endaligned \tag 2.12 $$
for $\sigma=\Re s>1/2$.  Let
$$l(s)=\sum_{d\in\Omega} \frac{h_d\ln \epsilon_d}{(d u_d^2)^s}.$$
By (2.11) and (2.12), we obtain that
$$l(s)-h(-i/2)\tag 2.13 $$
 is analytic in the half-plane $\Re\, s>1/2$.

Let $\bar L_N(s)$ be given as in Lemma 2.5.   Then
by (1.4), (4.4), (4.5), Theorem 4.3, Lemma 5.1, and
Lemma 5.3 of [7], we have that
$$\bar L_N(s)-h(-i/2) \tag 2.14$$
is an analytic function of $s$ in the half-plane
$\Re\, s> 1/2$ except for simple poles at
$s=1/2-i\kappa_j$, $j=1,2,\cdots, S$.  It follows from (2.13) and (2.14)
that
$$L(s)=\bar L_N(s)-l(s)$$
represents an analytic function of $s$ in the half-plane
$\Re\, s> 1/2$ except for simple poles at
$s_j=1/2-i\kappa_j$, $j=1,2,\cdots, S$.
Moreover, we have
$$m_j=2\,\text{Res}_{s=s_j} L(s)$$
for $j=1,2,\cdots, S$.

  This completes the proof of the theorem. \qed\enddemo

\demo{Proof of Corollary 2}  By Theorem 1 the series
$$\sum_{k|N}k\sum_{d\in \Omega,\, k|u_d}
\prod_{p^{2l}|(d, N/k)}p^l \prod_{p| N/k}\{1+\left({d\over p}\right)\}
\prod_{p|k}\{1-{1\over p}\left({d\over p}\right)\}
{h_d\ln \epsilon_d\over (du_d^2)^s} -\sum_{d\in\Omega}
\frac{h_d\ln \epsilon_d}{(d u_d^2)^s} $$
represents an analytic function of $s$ in the
half-plane $\Re s > 1/2$ except for having simple poles at
$s_j=\frac 12- i\kappa_j$, $j=1,2,\cdots, S$, where $(v_d, u_d)$ is
the smallest positive solution of Pell's equation $v^2-du^2=4$ and
the product on $p^{2l}$ is over all distinct primes $p$ with
$p^{2l}$ being the greatest even $p$-power factor of $(d,N/k)$.
 Moreover, we have
$$m_j=2\,\text{Res}_{s=s_j} L(s)$$
for $j=1,2,\cdots, S$.

Since $N$ is square free, we have
$$\prod_{p| N/k}\{1+\left({d\over p}\right)\}
\prod_{p|k}\{1-{1\over p}\left({d\over p}\right)\}
=\sum_{m|N} {\mu((m, k))\over (m,k)} \left({d\over m}\right) \tag 2.15$$
where $\mu(n)$ is the M\"obius function and $(m, k)$ denotes
the greatest common divisor of $m$ and $k$.
By using the identity (2.15), we can write
$$\aligned &\sum_{k|N}k\sum_{d\in \Omega,\, k|u_d}
\prod_{p^{2l}|(d, N/k)}p^l \prod_{p| N/k}\{1+\left({d\over p}\right)\}
\prod_{p|k}\{1-{1\over p}\left({d\over p}\right)\}
{h_d\ln \epsilon_d\over (du_d^2)^s} -\sum_{d\in\Omega}
\frac{h_d\ln \epsilon_d}{(d u_d^2)^s}\\
&= \sum_{m|N, k|N}k{\mu((m, k))\over (m,k)}\sum_{d\in \Omega,\, k|u_d}  \left({d\over m}\right)
{h_d\ln \epsilon_d\over (du_d^2)^s} -\sum_{d\in\Omega}
\frac{h_d\ln \epsilon_d}{(d u_d^2)^s}\\
&=\sum_{\underset{(m, k)\neq (1,1)}\to {m|N, \,k|N}}
k{\mu((m, k))\over (m,k)} \sum_{d\in \Omega,\,
k|u_d} \left({d\over m}\right) {h_d\ln \epsilon_d\over (du_d^2)^s}.
\endaligned  $$
It follows that then the series
$$L_1(s)=\sum_{\underset{(m, k)\neq (1,1)}\to {m|N, \,k|N}}
k{\mu((m, k))\over (m,k)} \sum_{d\in \Omega,\,
k|u_d} \left({d\over m}\right) {h_d\ln \epsilon_d\over (du_d^2)^s}$$
represents an analytic function of $s$ in the
half-plane $\Re s > 1/2$ except for having simple poles at those
points $s_j=\frac 12- i\kappa_j$, $j=1,2,\cdots, S$,
where $(v_d, u_d)$ is the smallest positive solution of Pell's
equation $v^2-du^2=4$.  Moreover, we have
$$m_j=2\,\text{Res}_{s=s_j} L_1(s)$$
for $j=1,2,\cdots, S$.

 This completes the proof of the corollary. \qed\enddemo

\enddocument